# A Reduced Forward Collatz Algorithm: How Binary Strings Change Their Length Under 3x+1

Richard Kaufman
rdkaufman01 at gmail dot com

**Abstract**

We developed an algorithm that easily goes from one odd number to the next odd number in binary representation for the reduced forward Collatz map (Syracuse function). The algorithm indicates when an odd number can grow or shrink to the next odd number based on the pattern of binary digits. The algorithm is also used to provide a simpler method for determining the change in binary string length for the reduced map than one found in the literature. Accordingly, an inspection of the binary digits for an odd number can determine the number of binary digits of the subsequent odd number. We also show some simple results for what the smallest number could be for a counterexample to the Collatz conjecture.

## 1. Introduction

The Collatz conjecture states that, for any positive integer $n$, the function

$$f(n) = \begin{cases} \dfrac{n}{2} & \text{if } n \text{ is even} \\ \dfrac{3n+1}{2} & \text{if } n \text{ is odd} \end{cases}$$

will always provide for a sequence of positive integers that reaches 1.

An alternative form of the function, called the reduced Collatz function, maps one odd number to the next odd number, such that only odd numbers are included in the Collatz sequence. In this paper, we present an algorithm, called Algorithm 1, that easily performs the reduced Collatz function in binary notation. The algorithm is then used to provide a method, called Method 1, for determining the change in binary string length for the reduced Collatz function. Method 1 appears to be simpler than one presented by Hew, which he calls Theorem 1.[1]

Hew (page 481) states that:
> A full description of how a bit string's length will change under application
> of the reduced Collatz map has yet to appear in the literature….
> This article provides a way to tell, by inspection, the change in length that a
> bit string will incur under the reduced Collatz map.

Theorem 1 provided by Hew (page 481) states:

> **Theorem 1.** Let $x$ be an odd positive integer and write the binary expansion of $x$ as $u[01]^n$ where $u$ is non-empty, lacks leading zeroes, and does not end with 01. Put $\hat{s} = (x + \frac{1}{3})/D(x)$ where $D(x) = 2^{\lfloor \log_2 x \rfloor + 1}$. Then $R(x)$ will be shorter than $x$ by the following number of bits:
>
> |  | u ends with 0 (so $n \geq 1$) | u ends with 1 (so $n \geq 0$) |
> |---|---|---|
> | $\hat{s} < 2/3$ | $2n-1$ bits | $2n$ bits |
> | $\hat{s} = 2/3$ | – | $2n$ bits |
> | $\hat{s} > 2/3$ | $2n-2$ bits | $2n-1$ bits |

Hew appears to be correct that his paper is the first to consider the change in binary length of the reduced Collatz function. However, the method presented in the present paper may be simpler to implement and inspect. In the course of developing this paper, we also observed another result for what the smallest number could be for a counterexample to the Collatz conjecture, which is also shown in the next section.

## 2. Main results

We propose a simple algorithm, called Algorithm 1, as an alternative for performing the reduced forward Collatz function. The reduced forward Collatz function goes from one odd number $d$ in the Collatz sequence to the next odd number $n$. The algorithm works on a binary representation of a number and can be verified by inspection.

**Algorithm 1.**
   Take odd $d$ in binary representation and put two "0"s in front of the leading "1", and then go through this binary string from right to left.
   First number to write is always 0.
   Set Flag = 1.
   If Flag = 0 then
   - If last two digits were "00", then write 0
     Elseif last two digits were "01" then write 1
     Elseif last two digits were "10" then write 1
     Elseif last two digits were "11" then write 0 and *set Flag = 1*
     End if
   Else If Flag = 1 then
   - If last two digits were "00", then write 1 and *set Flag = 0*
     Elseif last two digits were "01", then write 0
     Elseif last two digits were "10", then write 0
     Else if last two digits were "11", then write 1
     End If
   End If
   Stop when have exhausted all binary digits from right to left.

> Eliminate leading and trailing 0s for next odd number $n$.

The following example shows how Algorithm 1 is applied to the number 467.
Example A:

```
                467₁₀ =    1 1 1 0 1 0 0 1 1  Write number in binary
Number row                 0 0 1 1 1 0 1 0 0 1 1  Put 00 in front of binary. This is the "number row"
Steps:                                         0  First number to write is 0. Set Flag = 1.
                                             1 0  Flag = 1. Last two digits in number row were 11. So write 1.
                                           0 1 0  Flag = 1. Last two digits in number row were 01. So write 0.
                                         1 0 1 0  Flag = 1. Last two digits in number row were 00. So write 1. Set Flag = 0.
                                       1 1 0 1 0  Flag = 0. Last two digits in number row were 10. So write 1.
                                     1 1 1 0 1 0  Flag = 0. Last two digits in number row were 01. So write 1.
                                   1 1 1 1 0 1 0  Flag = 0. Last two digits in number row were 10. So write 1.
                                 0 1 1 1 1 0 1 0  Flag = 0. Last two digits in number row were 11. So write 0. Set Flag = 1.
                               1 0 1 1 1 1 0 1 0  Flag = 1. Last two digits in number row were 11. So write 1.
                             0 1 0 1 1 1 1 0 1 0  Flag = 1. Last two digits in number row were 01. So write 0.
                           1 0 1 0 1 1 1 1 0 1 0  Flag = 1. Last two digits in number row were 00. So write 1. Set Flag = 0. End
Next Number                1 0 1 0 1 1 1 1 0 1    Eliminate leading and trailing 0. This is 701 in base 10.
```

The following example shows each subsequent odd number of the reduced Collatz function starting with 31.
Example B:

| Sequence of odd Collatz numbers | Binary | Number of binary digits | Difference in number of binary digits from previous |
|---|---|---|---|
| 31 | 11111 | 5 |  |
| 47 | 101111 | 6 | 1 |
| 71 | 1000111 | 7 | 1 |
| 107 | 1101011 | 7 | 0 |
| 161 | 10100001 | 8 | 1 |
| 121 | 1111001 | 7 | -1 |
| 91 | 1011011 | 7 | 0 |
| 137 | 10001001 | 8 | 1 |
| 103 | 1100111 | 7 | -1 |
| 155 | 10011011 | 8 | 1 |
| 233 | 11101001 | 8 | 0 |
| 175 | 10101111 | 8 | 0 |
| 263 | 100000111 | 9 | 1 |
| 395 | 110001011 | 9 | 0 |
| 593 | 1001010001 | 10 | 1 |
| 445 | 110111101 | 9 | -1 |
| 167 | 10100111 | 8 | -1 |
| 251 | 11111011 | 8 | 0 |
| 377 | 101111001 | 9 | 1 |
| 283 | 100011011 | 9 | 0 |
| 425 | 110101001 | 9 | 0 |
| 319 | 100111111 | 9 | 0 |
| 479 | 111011111 | 9 | 0 |
| 719 | 1011001111 | 10 | 1 |
| 1079 | 10000110111 | 11 | 1 |
| 1619 | 11001010011 | 11 | 0 |
| 2429 | 100101111101 | 12 | 1 |
| 911 | 1110001111 | 10 | -2 |
| 1367 | 10101010111 | 11 | 1 |
| 2051 | 100000000011 | 12 | 1 |
| 3077 | 110000000101 | 12 | 0 |
| 577 | 1001000001 | 10 | -2 |
| 433 | 110110001 | 9 | -1 |
| 325 | 101000101 | 9 | 0 |
| 61 | 111101 | 6 | -3 |
| 23 | 10111 | 5 | -1 |
| 35 | 100011 | 6 | 1 |
| 53 | 110101 | 6 | 0 |
| 5 | 101 | 3 | -3 |
| 1 | 1 | 1 | -2 |

The following example shows each subsequent odd number starting with 63.
Example C:

| Base 10 | Binary (Start 63) | 12 | 11 | 10 | 9 | 8 | 7 | 6 | 5 | 4 | 3 | 2 | 1 |
|---|---|---|---|---|---|---|---|---|---|---|---|---|---|
| 63 | 111111 | | | | | | | 1 | 1 | 1 | 1 | 1 | 1 |
| 95 | 1011111 | | | | | | 1 | 0 | 1 | 1 | 1 | 1 | 1 |
| 143 | 10001111 | | | | | 1 | 0 | 0 | 0 | 1 | 1 | 1 | 1 |
| 215 | 11010111 | | | | | 1 | 1 | 0 | 1 | 0 | 1 | 1 | 1 |
| 323 | 101000011 | | | | 1 | 0 | 1 | 0 | 0 | 0 | 0 | 1 | 1 |
| 485 | 111100101 | | | | 1 | 1 | 1 | 1 | 0 | 0 | 1 | 0 | 1 |
| 91 | 1011011 | | | | | | 1 | 0 | 1 | 1 | 0 | 1 | 1 |
| 137 | 10001001 | | | | | 1 | 0 | 0 | 0 | 1 | 0 | 0 | 1 |
| 103 | 1100111 | | | | | | 1 | 1 | 0 | 0 | 1 | 1 | 1 |
| 155 | 10011011 | | | | | 1 | 0 | 0 | 1 | 1 | 0 | 1 | 1 |
| 233 | 11101001 | | | | | 1 | 1 | 1 | 0 | 1 | 0 | 0 | 1 |
| 175 | 10101111 | | | | | 1 | 0 | 1 | 0 | 1 | 1 | 1 | 1 |
| 263 | 100000111 | | | | 1 | 0 | 0 | 0 | 0 | 0 | 1 | 1 | 1 |
| 395 | 110001011 | | | | 1 | 1 | 0 | 0 | 0 | 1 | 0 | 1 | 1 |
| 593 | 1001010001 | | | 1 | 0 | 0 | 1 | 0 | 1 | 0 | 0 | 0 | 1 |
| 445 | 110111101 | | | | 1 | 1 | 0 | 1 | 1 | 1 | 1 | 0 | 1 |
| 167 | 10100111 | | | | | 1 | 0 | 1 | 0 | 0 | 1 | 1 | 1 |
| 251 | 11111011 | | | | | 1 | 1 | 1 | 1 | 1 | 0 | 1 | 1 |
| 377 | 101111001 | | | | 1 | 0 | 1 | 1 | 1 | 1 | 0 | 0 | 1 |
| 283 | 100011011 | | | | 1 | 0 | 0 | 0 | 1 | 1 | 0 | 1 | 1 |
| 425 | 110101001 | | | | 1 | 1 | 0 | 1 | 0 | 1 | 0 | 0 | 1 |
| 319 | 100111111 | | | | 1 | 0 | 0 | 1 | 1 | 1 | 1 | 1 | 1 |
| 479 | 111011111 | | | | 1 | 1 | 1 | 0 | 1 | 1 | 1 | 1 | 1 |
| 719 | 1011001111 | | | 1 | 0 | 1 | 1 | 0 | 0 | 1 | 1 | 1 | 1 |
| 1079 | 10000110111 | | 1 | 0 | 0 | 0 | 0 | 1 | 1 | 0 | 1 | 1 | 1 |
| 1619 | 11001010011 | | 1 | 1 | 0 | 0 | 1 | 0 | 1 | 0 | 0 | 1 | 1 |
| 2429 | 100101111101 | 1 | 0 | 0 | 1 | 0 | 1 | 1 | 1 | 1 | 1 | 0 | 1 |
| 911 | 1110001111 | | | 1 | 1 | 1 | 0 | 0 | 0 | 1 | 1 | 1 | 1 |
| 1367 | 10101010111 | | 1 | 0 | 1 | 0 | 1 | 0 | 1 | 0 | 1 | 1 | 1 |
| 2051 | 100000000011 | 1 | 0 | 0 | 0 | 0 | 0 | 0 | 0 | 0 | 0 | 1 | 1 |
| 3077 | 110000000101 | 1 | 1 | 0 | 0 | 0 | 0 | 0 | 0 | 0 | 1 | 0 | 1 |
| 577 | 1001000001 | | | 1 | 0 | 0 | 1 | 0 | 0 | 0 | 0 | 0 | 1 |
| 433 | 110110001 | | | | 1 | 1 | 0 | 1 | 1 | 0 | 0 | 0 | 1 |
| 325 | 101000101 | | | | 1 | 0 | 1 | 0 | 0 | 0 | 1 | 0 | 1 |
| 61 | 111101 | | | | | | | 1 | 1 | 1 | 1 | 0 | 1 |
| 23 | 10111 | | | | | | | | 1 | 0 | 1 | 1 | 1 |
| 35 | 100011 | | | | | | | 1 | 0 | 0 | 0 | 1 | 1 |
| 53 | 110101 | | | | | | | 1 | 1 | 0 | 1 | 0 | 1 |
| 5 | 101 | | | | | | | | | | 1 | 0 | 1 |
| 1 | 1 | | | | | | | | | | | | 1 |

The colored cells show some diagonal patterns that were also observed for all Collatz sequences:
- Yellow shows the first 0 (read right to left)
- Blue shows adjacent alternating bits on the diagonal
- Red shows adjacent bits of either 0 or 1 until each diagonal terminates.

The algorithm indicates when an odd number can grow or shrink to the next odd number based on the pattern of binary digits. The following are simple observations based upon Algorithm 1:
1. When reading binary representation from *left to right*:
    a. Leading binary digit 1 can only move to the left by one place from one odd number $d$ in the sequence to the next odd number $n$. This can occur when $d$ has "11" in it without "00" first, starting from the leading digit 1.

b. Leading binary digit 1 can move to the right by more than one place from one number $d$ to the next odd number $n$ (this will be specified later in Method 1).
c. When an odd number ends in "01", then the next odd number can have the same number of digits or less.
2. When reading binary representation from *right to left*:
a. The first binary digit 0 moves one place to the right for each subsequential odd number, until "01" is reached.

This algorithm provides more insight then just using binary addition[2].

We can also skip numbers in the Collatz sequence by going from one odd number to the odd number that follows the *next* even number. To do so, we will use Theorem A, which is discussed next.

**Theorem A**: Every positive odd number can be expressed as $n = g2^k - 1$, for some positive odd number $g$ and some positive integer $k$. For an odd number in the forward Collatz sequence (not the reduced forward Collatz sequence with only odd numbers), the *highest subsequent number* (which is even) *in a continually ascending* forward Collatz sequence will be $g3^k - 1$.

Before we prove this, we provide some examples:

*Example* 1: For odd $n = 79 = 5(2)^4 - 1$, the highest subsequent number in a continually ascending forward Collatz sequence is $5(3)^4 - 1 = 404$.

*Example* 2: For odd $n = 127 = (2)^7 - 1$, the highest subsequent number in a continually ascending forward Collatz sequence is $(3)^7 - 1 = 2186$.

The following example shows how this can reduce the number of calculations starting with 63 to reach 1.
Example D:

| Base 10 | Binary (Start 63) | 12 | 11 | 10 | 9 | 8 | 7 | 6 | 5 | 4 | 3 | 2 | 1 |
|---|---|---|---|---|---|---|---|---|---|---|---|---|---|
| 63 | 111111 |  |  |  |  |  |  | 1 | 1 | 1 | 1 | 1 | 1 |
| 91 | 1011011 |  |  |  |  |  | 1 | 0 | 1 | 1 | 0 | 1 | 1 |
| 103 | 1100111 |  |  |  |  |  | 1 | 1 | 0 | 0 | 1 | 1 | 1 |
| 175 | 10101111 |  |  |  |  | 1 | 0 | 1 | 0 | 1 | 1 | 1 | 1 |
| 445 | 110111101 |  |  |  | 1 | 1 | 0 | 1 | 1 | 1 | 1 | 0 | 1 |
| 167 | 10100111 |  |  |  |  | 1 | 0 | 1 | 0 | 0 | 1 | 1 | 1 |
| 283 | 100011011 |  |  |  | 1 | 0 | 0 | 0 | 1 | 1 | 0 | 1 | 1 |
| 319 | 100111111 |  |  |  | 1 | 0 | 0 | 1 | 1 | 1 | 1 | 1 | 1 |
| 911 | 1110001111 |  |  | 1 | 1 | 1 | 0 | 0 | 0 | 1 | 1 | 1 | 1 |
| 577 | 1001000001 |  |  | 1 | 0 | 0 | 1 | 0 | 0 | 0 | 0 | 0 | 1 |
| 433 | 110110001 |  |  |  | 1 | 1 | 0 | 1 | 1 | 0 | 0 | 0 | 1 |
| 325 | 101000101 |  |  |  | 1 | 0 | 1 | 0 | 0 | 0 | 1 | 0 | 1 |
| 61 | 111101 |  |  |  |  |  |  | 1 | 1 | 1 | 1 | 0 | 1 |
| 23 | 10111 |  |  |  |  |  |  |  | 1 | 0 | 1 | 1 | 1 |
| 5 | 101 |  |  |  |  |  |  |  |  |  | 1 | 0 | 1 |
| 1 | 1 |  |  |  |  |  |  |  |  |  |  |  | 1 |

This only requires 16 rows instead of the 40 used previously. This is based on Theorem A and provides additional insight into the numbers used in the parity sequence.[3]

**Proof of Theorem A by induction**:
*Base case*. For $k = 1$, then $n = g2^1 - 1$ which is odd. The forward Collatz function is $f(n) = \frac{3(g2^1-1)+1}{2} = g3^1 - 1$. This number is even, since $g$, $g3^k$, and 1 are all odd. So, the next number in the forward Collatz sequence would be smaller. Therefore, the base case is true.
*Inductive step*. Assume $n = g2^k - 1$ for some positive integer $k$ and odd number $g$, such that the highest subsequent number in a continually ascending forward Collatz sequence is $g3^k - 1$. Show that this is also true for $k + 1$.
By the inductive assumption, $f(g2^k - 1) = g3^k - 1$. So $f(g2^{k+1} - 1) = f((2g)2^k - 1) = (2g)3^k - 1$, which is an odd number equal to $(3^k g)2^1 - 1$. Again, by the inductive assumption, $f\left((3^k g)2^1 - 1\right) = (3^k g)3^1 - 1 = g3^{k+1} - 1$. Since $g$ is odd, then we have the desired result of $g3^{k+1} - 1$, which is even, so the next number would be smaller. Q.E.D.

By inspection (*i.e.*, see Example C and D), when odd $n$ is written in binary form, it is apparent that $k$ is the number of continuous 1s, from right to left (*i.e.*, it takes $k$ steps in the reduced forward Collatz sequence for the first 0 bit, from right to left, to move all the way to the right). So $k > 0$.

To go from odd number $n_i$ to the odd number $n_{i+1}$ following the *next* even number in the sequence, then $n_{i+1} = g3^k - 1$, where all powers of 2 have been removed.
Since $n_i + 1 = g2^k$, then: $n_{i+1} = g3^k - 1 = \frac{(g2^k)(g3^k)}{(g2^k)} - 1 = (n_i + 1)\left(\frac{3^k}{2^k}\right) - 1 = (n_i + 1)\left(\frac{3}{2}\right)^k - 1.$

Therefore, the following Equation 1 takes one odd number $n_i$ to the odd number $n_{i+1}$ following the *next* even number in the sequence, then:

$$n_{i+1} = (n_i + 1)\left(\frac{3}{2}\right)^k - 1, \text{ where all powers of 2 have been removed.}$$

(1)

Since $n_{i+1}$ is an odd number where all powers of 2 have been removed, then $w$ is the number of continuous 1s in $(n_i + 1)\left(\frac{3}{2}\right)^k - 1$, from right to left, when written in binary. Or,

$$n_{i+1} = \frac{(n_i + 1)\left(\frac{3}{2}\right)^k - 1}{2^w}$$

(2)

This can be rewritten as:

$$\frac{2^w n_{i+1} + 1}{n_i + 1} = \left(\frac{3}{2}\right)^k$$

where $k > 0$ and $w > 0$.

(3)

Hew provided Theorem 1, referenced earlier, to determine the change in binary string length for a reduced Collatz map. Here we use Algorithm 1 to develop a simpler method, called Method 1 to distinguish is from Hew's Theorem 1.

**Method 1.**
Determination of the number of binary digits from one odd number $d$ to the next odd number $n$ in the Collatz sequence:
- Take odd $d$ in binary representation and put **one** "0" in front of the leading "1", and then go through this binary string from right to left.
- Set Flag = 1
- If Flag = 0 and encounter "11" then set Flag = 1
  Elseif Flag = 1 and encounter "00" then set Flag = 0
  End if
- Stop when have exhausted all binary digits from right to left.
- Num = Number of contiguous alternating digits starting from the right end of $d$ (note first digit 1 always counts as 1).

$n$ will have (1+ Flag – Num) number of digits more than $d$.

Method 1 is easy to implement. It shows that an inspection of the digits of an odd number can determine the change in the number of digits to the next odd number in the reduced Collatz sequence.

**3. Additional results**

In the course of developing this paper, we also observed some other results for the Collatz conjecture, which are shown in this section.

**Theorem S**: If the Collatz conjecture is false, then the sequence has a smallest positive integer $n$. Then $n$ cannot be of any of the following forms:
- $2w$
- $1 + 4w$
- $2 + 3w$
- $3 + 16w$

for any whole number $w$. For such an $n$, then the first bullet shows that $n$ would be odd.

Let us consider a proof by contradiction for the Collatz conjecture.
Without loss of generality, we will assume that all positive integers *below* integer $n$ reach 1. But, that positive integer $n$ does not reach 1. [In other words, $n$ is the first positive integer that does not reach 1.] We will aim to show that a contradiction would result, which would prove the conjecture.

By the assumption, the (forward or backwards) Collatz sequence of integers from $n$ cannot have an integer less than $n$ or it would reach 1 (contradicting the initial assumption). Therefore, $n$ is the smallest integer in its Collatz sequence.

Now, if $n$ were even, then $n/2$ would be the next number in the sequence and smaller, which is a contradiction of the initial assumption. So, $n$ must be odd.

Restriction 1. $n$ is odd. So, $n$ is *not* of the form $2w$, for whole number $w$.

Since $n$ is the smallest number in its forward or backward sequence, then the number in the sequence prior to $n$ must have been larger. If the prior number $m$ from the Collatz function was odd (*i.e.,* $m$ such that $f(m) = (3m + 1)/2 = n$), then $m = (2n - 1)/3$, which is smaller than $n$, and impossible. So, $m$ must be even. Since $m$ cannot be odd, then from the Collatz function $m = (2n - 1)/3 \neq 1 + 2w$, where $w$ is a whole number. This means that $n \neq ((3 + 6w) + 1)/2$, or $n \neq 2 + 3w$. Therefore, $n$ cannot be a number such as 2, 5, 8, 11, 14,…, [numbers of the form $2 + 3w$ for whole number $w$]:

Restriction 2. $n$ is *not* of the form $2 + 3w$, for whole number $w$.

Putting odd $n$ into the Collatz function, then $o = \frac{3n+1}{2}$ is the next positive integer in the sequence. Now if $o$ was even, the subsequent number would be $((3n + 1)/2)/2 = (3n + 1)/4$, which would be less than $n$, which is impossible by the initial assumption.

So, $o$ must be odd. Since $o$ cannot be even, then $o = (3n + 1)/2 \neq 2w$, where $w$ is a positive integer. This means that $n \neq (4w - 1)/3$. Therefore, $n$ cannot be a number such as 1, 5, 9, 13, 17,… [numbers of the form $1 + 4w$ for whole number $w$]:

*Restriction 3*: $n$ is *not* of the form $1 + 4w$, for whole number $w$.

Putting odd $o$ into the Collatz function, then $p = \frac{3 \cdot o + 1}{2} = \frac{\left(3\left(\frac{3n+1}{2}\right)+1\right)}{2} = (9n+5)/4$ is the next positive integer in the sequence. Now $p$ cannot be divisible by $2^2$ or the sequence will result in a number less than $n$. Therefore, positive integer $p = (9n+5)/4 \neq 4w$, where $w$ is a whole number. This means that $n \neq (16w - 5)/9$. Therefore, $n$ cannot be a number such as 3, 19, 35, 51, 67, … [numbers of the form $3 + 16w$ for whole number $w$].

*Restriction 4*: $n$ is *not* of the form $3 + 16w$, for whole number $w$.

Combining the previous restrictions, then $n$ is *not* of the form:
$2w, 1 + 4w, 2 + 3w, 3 + 16w$
for any whole number $w$.

With these restrictions, then $n$ could **not** be an odd number such as 3, 5, 9, 11, 13, 17, 19, 21, 23, 25, 29, 33, 35, 37, 41,…

Alternatively, using this approach, then $n$ could be an odd number such as 7, 15, 27, 31, 39, 43, 55, 63, 75, 79, 87, 91, 95, 103, 107, 111, 117, 119,… Of course, this means that no forward Collatz sequence starting from any of these numbers could have any odd number smaller than the initial number $n$, since this would eventually reach 1 (by the initial assumption that $n$ is a smallest number that violates the Collatz conjecture).

## 4. Conclusion

We present a simple algorithm for performing the forward Collatz function and use it to develop a simpler method for determining the change in binary string length than one presented by Hew. Hew stated (page 486) that, "The reformulation yields insights that are useful and beautiful" and also "speculates that more results will come from viewing the Collatz function" in this way. We agree.

## Acknowledgement
The author thanks Pablo Castañeda for corrections to this paper. He has agreed to collaborate on a new paper using these results.

**Funding details:** There was no funding provided for this paper.
**Disclosure statement:** The authors reports there are no competing interests to declare.

---

[1] Patrick Chisan Hew (2021) Collatz on the Dyadic Rationals in [0.5, 1) with Fractals: How Bit Strings Change Their Length Under 3x + 1, Experimental Mathematics, 30:4, 481-488, DOI: 10.1080/10586458.2019.1577765

---

[2] "Collatz Conjecture." *As an abstract machine that computes in base two*. Wikipedia, 9 Sept. 2020, en.wikipedia.org/wiki/Collatz_conjecture. (accessed 6/19/2022).

[3] "Collatz Conjecture." *As a parity sequence*. Wikipedia, 9 Sept. 2020, en.wikipedia.org/wiki/Collatz_conjecture. (accessed 6/19/2022).